\documentclass[11pt, a4paper, reqno]{amsart}
\usepackage{latexsym}
\usepackage{tcolorbox}

\usepackage{amssymb, amsmath, amsthm, amsrefs, xcolor}
\usepackage{mathtools}
\usepackage{hyperref}
\usepackage{mathrsfs}
\usepackage[matha, mathx]{mathabx}
\usepackage[shortlabels]{enumitem}
\setlist{leftmargin=*,topsep=0.2em}

\usepackage{hyperref} 
\hypersetup{
    colorlinks=false,       
    linkcolor=blue,          
    citecolor=magenta,        
    filecolor=magenta,      
    urlcolor=cyan           
}

\newtheorem{thm}{Theorem}

\newcommand{\dd}{\mathrm{d}}

\begin{document}
\title{A note on bilinear multipliers with convex singularities}
\author{ Valentina Ciccone }
\address{Institute of Mathematics, Polish Academy of Sciences, Śniadeckich 8, 00-656 Warszawa, Poland}
\email{vciccone@impan.pl}

\begin{abstract}
    We study bounds in the local $L^2$ range of exponents for bilinear multipliers whose symbol is the characteristic function of the epigraph of certain convex curves.
    We realize these bounds as a consequence of estimates that we establish, via simple arguments, for the associated exotic paraproducts. 
    As a further application, we observe bounds beyond the local $L^2$ range for bilinear multipliers whose symbol is the characteristic function of the epigraph of convex polygonal curves associated with these paraproducts.
\end{abstract}

\thanks{The author was supported by the Deutsche Forschungsgemeinschaft (DFG, German Research Foundation) under
Germany’s Excellence Strategy – EXC2047/1 390685813 as well as SFB 1060, and by the National Science Centre, Poland, grant Sonata Bis 2022/46/E/ST1/00036.} 

\maketitle

\section{Introduction}
Let $m$ be a bounded measurable function on $\mathbb{R}^2$. For a pair of Schwartz functions $f,g\in\mathscr{S}(\mathbb{R})$ we define the bilinear multiplier operator
\begin{align}\label{multiplier_operator}
    B_m(f,g)(x):=\int \int_{\mathbb{R}^2} m(\xi,\eta)\widehat{f}(\xi)\widehat{g}(\eta) e^{2\pi i (\xi+\eta)x} \dd\xi \dd\eta~.
\end{align}
We are interested in bounds for $B_m$ of the form
\begin{align}\label{wanted_bound}
    \Vert B_m(f,g)\Vert_{L^{p_3'}(\mathbb{R})} \leq C(p_1,p_2,m)\Vert f\Vert_{L^{p_1}(\mathbb{R})} \Vert g \Vert_{L^{p_2}(\mathbb{R})}~, 
\end{align}
mostly in the local $L^2$ region
\begin{align}\label{local_L2}
    2 \leq  p_1, p_2, p_3 <\infty~, \quad  \frac{1}{p_1} +\frac{1}{p_2} +\frac{1}{p_3}=1~,
\end{align}
where we use the notation $p'_3:=p_3/(p_3-1)$ for the dual exponent.

In this note, the Fourier multiplier symbol $m$ will always be the indicator function of some set $\Gamma\subseteq \mathbb{R}^2$, $m=\mathbf{1}_\Gamma$. If $\Gamma=\mathbb{R}^2$ we immediately recover the pointwise product operator mapping $(f,g)$ in $fg$. On the other hand, choosing $\Gamma$ to be a half-plane yields the bilinear Hilbert transform. Bounds of the form \eqref{wanted_bound}
in the open local $L^2$ range for the bilinear Hilbert transform have been first established by Lacey and Thiele in \cite{LT97,LT97a}. These results have been later extended beyond the local $L^2$ range in \cite{LT98,LT99,DPT16} and uniform estimates with respect to the slope of the associated half-plane have been established in \cite{T02,GL04,L06,UW22}.
The case of $\Gamma=\mathbb{D}$, the unit disc, has been studied by Grafakos and Li in \cite{GL06} where they established the boundedness of the associated bilinear multiplier operator in the local $L^2$ region \eqref{local_L2}. Due to the invariance under translation and dilation of bilinear multiplier norms, the local $L^2$ boundedness of the bilinear disc multiplier implies uniform bounds for the bilinear Hilbert transform in the same range of exponents \cite{GL06}.
Bounds for the bilinear multiplier associated with the epigraph of a parabola have been proved by Muscalu in \cite{M00}. 
Moreover, bounds for bilinear multipliers associated with smooth curves with bounded slopes were also established in \cite{M00}.
The boundedness of the bilinear multiplier associated with a lacunary polygon inscribed in the unit disc has been proved by Demeter and Gautam in \cite{DG12}.
More recently, Saari and Thiele \cite{ST23} have established boundedness, in the open local $L^2$ range, for bilinear multiplier operators whose symbol is the characteristic function of the epigraph of an exponential curve. The core of the analysis in \cite{ST23} consists in studying the boundedness in the open local $L^2$ range of an ad hoc paraproduct associated with a staircase set inscribed in the exponential curve. The passage from estimates for the paraproduct to estimates for the multiplier has then been achieved by invoking orthogonality-type arguments and analysis of uniform boundedness for bilinear multipliers associated with certain segments of curves with constrained slope from \cite{GL06}, see also \cite{M00}. The decomposition of bilinear multipliers into paraproducts is a well-established technique. Besides the aforementioned works, we also mention \cite{L08,MTT02}, and for a more exhaustive list of classical references, we refer to the textbook \cite{MS13}.
In general, it is an open problem to find a general description of sets in $\mathbb{R}^2$ whose characteristic function defines a bounded bilinear multiplier operator \cite{GL06}, begging the question of whether convexity of the set is, or not, a sufficient condition \cite{ST23}. For example, the answer to this latter question is known to be negative for bilinear multipliers in $\mathbb{R}^{2d}$, $d\geq 2$, outside the local $L^2$ region  \cite{DG07,GR10,G12}. 

In this note, we study the boundedness -- primarily in the local $L^2$ region of exponents -- of bilinear multipliers whose symbol is the characteristic function of the epigraph of certain further convex curves.
In doing so, we also recover, via a unified approach, the boundedness of bilinear multipliers associated with 
some previously studied convex curves. 
Inspired by the approach in \cite{ST23}, at the heart of our analysis there are bounds that we establish for certain families of exotic paraproducts associated with staircase sets that we will introduce shortly. We stress that we establish such bounds via extremely simple arguments and that such bounds hold naturally beyond the local $L^2$ range.

Model case families of exotic paraproducts that we consider can be constructed as follows.
Let $\lbrace a_j\rbrace_{j\in\mathbb{N}_0}, \; \lbrace b_j \rbrace_{j\in\mathbb{N}_0}$ be two strictly monotone sequences of real numbers. Sometimes we will use the shorthand notation ${\mathtt{a}}:=\lbrace a_j\rbrace_{j\in\mathbb{N}_0}, \; \mathtt{b}:=\lbrace b_j \rbrace_{j\in\mathbb{N}_0}$. We may assume that these sequences are strictly decreasing. 
We define the multiplier symbol $ m_{\mathtt{a,b}}$ associated with two strictly decreasing sequences $\lbrace a_j\rbrace_{j\in\mathbb{N}_0}, \, \lbrace b_j\rbrace_{j\in\mathbb{N}_0}$ as
    \begin{align}\label{defi_m_a_b}
    m_{\mathtt{a,b}}(\xi,\eta) := \sum_{j\in\mathbb{N}} \mathbf{1}_{[a_{j+1},a_j)}(\xi)\mathbf{1}_{[b_j,b_0)}(\eta)
    \end{align}
    and the associated bilinear multiplier operator (paraproduct) as
    \begin{align}\label{defi_B_m_a_b}
        B_{m_{\mathtt{a,b}}}(f,g)(x):=\int\int_{\mathbb{R}^2} m_{\mathtt{a,b}}(\xi,\eta)\widehat{f}(\xi)\widehat{g}(\eta)e^{2\pi i (\xi+\eta)x} \dd \xi \dd \eta~.
    \end{align}
Given two intervals $A,B\subseteq\mathbb{R}$ we denote by $-A-B$ the Minkowski sum $-A-B :=\lbrace -a-b:\, a\in A, \, b\in B\,  \rbrace$.
With this notation, we introduce our key assumptions on the sequences $\lbrace a_j\rbrace_{j\in\mathbb{N}_0}, \, \lbrace b_j\rbrace_{j\in\mathbb{N}_0}$. We will assume that at least one of the following holds:

\begin{itemize}[leftmargin=4em,topsep=7pt]
    \item[(Hyp 1)] The collection $\mathcal{I}_{\mathtt{a,b}}:=\lbrace -[a_{j+1}, a_j) - [b_{j},b_0) \,\rbrace_{j\in\mathbb{N}}$ can be split into a finite number $n(\mathtt{a,b})$ of subcollections 
whose intervals are pairwise disjoint.
\item[(Hyp 2)] The sequence $\lbrace b_j \rbrace_{j\in\mathbb{N}_0}$ is such that $\lim_{j\rightarrow\infty}b_j=b_\infty$, for some $b_\infty\in\mathbb{R}$, and the collection of intervals $\mathcal{J}_{\mathtt{a,b}}:=\lbrace -[a_{j+1},a_j)-(b_\infty,b_j) \, \rbrace_{j\in\mathbb{N}_0}$ can be split into a finite number of subcollections $n(\mathtt{a,b})$ whose intervals are pairwise disjoint.
\end{itemize}

To state our boundedness result for the paraproducts we need to recall a couple of definitions. 
We say that a closed null set $E\subset \mathbb{R}$ is Littlewood--Paley if for all $1<p<\infty$ there exist $c_{p,E},\, C_{p,E} > 0$ such that for all $f\in L^p(\mathbb{R})$ it holds that
 $$c_{p,E} \Vert f \Vert_{L^p} \leq \Vert \lbrace (\mathbf{1}_\omega\widehat{f})^\vee \rbrace_{\omega\in\Omega}\Vert_{L^p(\ell^2)} \leq C_{p,E} \Vert f \Vert_{L^p}$$
 where $\Omega$ is the countable collection of connected components of $\mathbb{R}\setminus E$, see e.g. \cite{HK89}. Moreover, for the purpose of this note, we say that  a collection of disjoint intervals $\Theta$ is Littlewood--Paley if for all $1<p<\infty$ there exist $C_{p,\Theta} > 0$ such that for all $f\in L^p(\mathbb{R})$ it holds that
    $$\Vert \lbrace (\mathbf{1}_\theta \widehat{f})^\vee \rbrace_{\theta\in\Theta} \Vert_{L^p(\ell^2)} \leq C_{p,\Theta} \Vert f \Vert_{L^p}~.$$
\begin{thm}\label{thm_general_paraproduct}
    Let $\lbrace a_j\rbrace_{j\in\mathbb{N}_0}, \;\lbrace b_j \rbrace_{j\in\mathbb{N}_0}$ be two strictly decreasing sequences of real numbers satisfying hypothesis $\mathrm{(Hyp \; 1)}$ -- respectively, $\mathrm{(Hyp \; 2)}$. Let $B_{m_{\mathtt{a,b}}}$ be the associated multiplier operator defined as in \eqref{defi_B_m_a_b}. 
    Then, the a priori estimate
    \begin{equation}\label{main_estimate_general_paraproduct}
        \Vert B_{m_{\mathtt{a,b}}}(f,g) \Vert_{p_3'} \leq C(p_1,p_2,\mathtt{a,b}) \Vert f \Vert_{p_1} \Vert g \Vert_{p_2}~,
    \end{equation}
    holds for all Schwartz functions $f,g\in\mathscr{S}(\mathbb{R})$ whenever
    \begin{equation*}
        2\leq p_1, p_3 < \infty~, \quad 1<p_2<\infty~, \quad \text{and} \quad \frac{1}{p_1}+\frac{1}{p_2}+\frac{1}{p_3} =1~,
    \end{equation*}
    and, in such range of exponents, $C(p_1,p_2,\mathtt{a,b})=C(p_1,p_2,n(\mathtt{a,b}))$.
    Moreover, if the sequence $\lbrace a_j\rbrace_{j\in\mathbb{N}_0},$ is Littlewood--Paley then \eqref{main_estimate_general_paraproduct} holds whenever
    \begin{equation*}
        1 < p_1, p_2 < \infty~, \quad 2\leq p_3 <\infty~, \quad \text{and} \quad \frac{1}{p_1}+\frac{1}{p_2}+\frac{1}{p_3} =1~.
    \end{equation*}
    If, in addition, the collection of intervals $\mathcal{I}_{\mathtt{a,b}}$ (respectively, $\mathcal{J}_{\mathtt{a,b}}$) can be split into a finite number of subcollections which are also Littlewood--Paley then \eqref{main_estimate_general_paraproduct} holds whenever 
    \begin{equation*}
        1 < p_1, p_2, p_3 < \infty~,\quad \text{and} \quad \frac{1}{p_1}+\frac{1}{p_2}+\frac{1}{p_3} =1~.
    \end{equation*}
\end{thm}
The result can be naturally adapted to handle the case of paraproducts associated with strictly increasing sequences, as well as the case of a strictly increasing and a strictly decreasing sequence. Moreover, the asymmetry in the range of the exponents can be reversed by 
inverting the role of $\xi$ and $\eta$. We will give an example of this in Section \ref{section_proof_paraproduct}.
We observe that to obtain the full range $1<p_2<\infty$ no Littlewood--Paley type assumption is needed for the sequence $\lbrace b_j \rbrace_{j\in\mathbb{N}_0}$. This is a byproduct of our key assumptions (Hyp 1), (Hyp 2). 
These hypotheses will be effective in dealing with paraproducts associated with staircase sets inscribed by certain sufficiently flat (or sufficiently steep, by reversing the role of $\xi$ and $\eta$) portions of convex curves. For example, the staircase paraproduct considered in \cite{ST23}, which leads to the boundedness of bilinear multipliers associated with exponential curves, falls into the class of paraproducts covered by Theorem \ref{thm_general_paraproduct}. To wit, Theorem \ref{thm_general_paraproduct} provides, as a particular case, a different proof of \cite[Theorem~1.1]{ST23}, see the first example in Section \ref{section_examples}.
Moreover, staircase paraproducts satisfying hypothesis (Hyp 2) will be well suited to study bilinear multipliers associated with segments of certain convex curves in the proximity of critical points.

We present two model case applications of Theorem \ref{thm_general_paraproduct}. The first one implies, as a particular case, the boundedness of the bilinear multiplier associated with the rectangular hyperbola $\xi\eta=1$, $\xi>0$.
\begin{thm}\label{thm_1_over_x_power}
    Consider the curve $\gamma(\xi):= |\xi|^{-c}$ for $\xi<0$, $c>0$, and let $m$ be the characteristic function of the set
    $$\lbrace (\xi,\eta): \,  \xi < 0, \, \eta \geq \gamma(\xi) \rbrace~.$$
    Then the bilinear multiplier operator $B_m$ defined as in \eqref{multiplier_operator}  satisfies the a priori bound \eqref{wanted_bound} for all triple of exponents $p_1,p_2,p_3$ in the local $L^2$ range \eqref{local_L2}.
\end{thm}
The second example implies the boundedness of the bilinear multiplier associated with the conjugate hyperbola $\xi^2-\eta^2=-1$.
\begin{thm}\label{thm_power_positive}
    Consider the curve $\gamma(\xi):= \sqrt{1+\xi^2}$, on $(-\tfrac{1}{\sqrt{3}},\tfrac{1}{\sqrt{3}})$ and 
    let $m$ be the characteristic function of the set
    $$\lbrace (\xi,\eta): \,  \xi\in (-\tfrac{1}{\sqrt{3}},\tfrac{1}{\sqrt{3}}), \, \eta \geq \gamma(\xi) \rbrace~.$$
    Then the bilinear multiplier operator $B_m$ defined as in \eqref{multiplier_operator} 
    satisfies the a priori bound \eqref{wanted_bound} for all triple of exponents $p_1,p_2,p_3$ in the local $L^2$ range \eqref{local_L2}.
\end{thm}
 We stress that the curves considered in Theorem \ref{thm_1_over_x_power} and Theorem \ref{thm_power_positive} are of interest to us because of their behavior near the so-called degenerate directions. However, the arguments in the proofs depend solely on the convexity of the functions and on certain properties of the sequences $\lbrace (\gamma ')^{-1}(2^{-j}) \rbrace_j$ (which will become clear shortly), any other curve with similar properties being equally interesting to us.
In Section \ref{section_examples} we collect more examples 
of sets for which the boundedness of the associated multiplier operator can be addressed, in a unified manner, via the same strategy.
In particular, we provide further examples and descriptions of convex and more general curves that circumscribe a staircase set whose associated paraproduct is covered by Theorem \ref{thm_general_paraproduct} and its variants and for which the passage from bounds for the paraproduct to bounds for the associated multiplier in the local $L^2$ range of exponents can be
completed as in \cite{GL06,M00,ST23}. Our examples include:
\begin{itemize}[leftmargin=2em,topsep=7pt]
    \item[(i)] any strictly convex, increasing, curve $\gamma\in\mathscr{C}^1$, with $0<\gamma'<1$, for which the sequence $\lbrace (\gamma ')^{-1}(2^{-j}) \rbrace_j$ is a convex sequence, see definition in Section \ref{section_examples};
    \item[(ii)] any strictly convex, increasing, curve $\gamma\in\mathscr{C}^1$, with $0<\gamma'<1$, for which the sequence $\lbrace (\gamma ')^{-1}(2^{-j}) \rbrace_j$ is a lacunary sequence, see definition in Section \ref{section_examples};
    \item[(iii)]  any strictly convex, increasing curve $\gamma\in\mathscr{C}^1$ on $(a,b)$, such that $\lim_{\xi\rightarrow a}$ $\gamma(\xi)=b_\infty$, for some $b_\infty \in\mathbb{R}$, with $0<\gamma'<1$, and for which  the sequence $\lbrace (\gamma ')^{-1}(2^{-j}) \rbrace_j$ is a concave sequence, see definition in Section \ref{section_examples};
    \item[(iv)] any, non necessarily convex, polygonal curve circumscribing a staircase paraproduct satisfying (Hyp 2). 
\end{itemize} 
As usual, the cases in the above list are understood up to translations, dilations, and finite decompositions of the multiplier symbols, as well as up to reversing the role of $\xi$ and $\eta$ and up to symmetries in the arguments.
A simple instance of a curve addressed by point (i) is the exponential curve $\gamma(\xi)=2^\xi$. The curves considered in Theorem \ref{thm_1_over_x_power} are examples of curves covered by point (ii). The curve $\gamma(\xi)=\xi^c/c$, $c>1$, on $(0,1)$, the curve $\gamma(\xi)=-\sqrt{1-\xi^2}$ on $(0,\tfrac{1}{\sqrt{2}})$, and the curve considered in Theorem \ref{thm_power_positive} are examples of curves that fall into the description in point (iii).

Finally, we record the following application to polygonal curves which extends beyond the local $L^2$ range. In this case, the passage from the paraproduct to the bilinear multiplier associated with the piecewise linear, convex curve follows closely the lines of the work of Demeter and Gautam \cite{DG12} where the boundedness beyond the local $L^2$ range of exponents has been established for the lacunary polygon multiplier.
\begin{thm}\label{thm_polygonal}
Let $\lbrace a_j\rbrace_{j\in\mathbb{N}_0}$, $\lbrace b_j \rbrace_{j\in\mathbb{N}_0}$ be two sequences as those considered in Theorem \ref{thm_general_paraproduct} and satisfying $\mathrm{(Hyp\; 2)}$. Moreover, assume that $\lbrace a_j\rbrace_{j\in\mathbb{N}_0}$, $\lbrace b_j \rbrace_{j\in\mathbb{N}_0}$ are such that the piecewise linear curve $\gamma$ of vertexes $\lbrace (a_j,b_j) \rbrace_{j\in\mathbb{N}_0}$ is convex and each line segment has slope in $(0,1)$.
Let $m$ be the characteristic function of the epigraph of $\gamma$.
Then, the corresponding bilinear multiplier operator $B_m$ defined as in \eqref{multiplier_operator} satisfies the a priori bound \eqref{wanted_bound} for all triples of exponents $p_1,p_2,p_3$ such that
$$2\leq p_1, p_3 < \infty~, \quad 1<p_2<\infty~, \quad \text{and} \quad \frac{1}{p_1}+\frac{1}{p_2}+\frac{1}{p_3} =1~.$$
\end{thm}
The passage from bounds for the paraproduct to bounds for the multiplier is naturally linked with uniform bounds for the bilinear Hilbert transform. This passage, beyond the local $L^2$ range of exponents and for the case of more general polygonal curves and convex curves, will be addressed elsewhere.

We conclude the introduction by mentioning some additional related works. Bilinear multipliers with non-smooth symbol and singularities along lines have been studied by Gilbert and Nahmod in \cite{GN00}. More recently, a sharp Hörmander condition for bilinear multipliers with singularities on a finite union of Lipschitz curves with tangencies away from degenerate directions has been established by Chen, Hsu, and Lin in \cite{CHL24}.

\subsection*{Structure} The paper is organized as follows. In Section \ref{section_proof_paraproduct} we prove Theorem \ref{thm_general_paraproduct} and we state a simple variant of it. In Section \ref{section_proof_1_over_x_power} we prove Theorem \ref{thm_1_over_x_power} and Theorem \ref{thm_power_positive}. Section \ref{section_examples} contains the aforementioned additional examples and generalizations.  Section \ref{section_polygonal} contains the proof of Theorem \ref{thm_polygonal}.

\section{Boundedness of the staircase paraproducts}\label{section_proof_paraproduct}
In this Section we prove Theorem \ref{thm_general_paraproduct}  and we give an example of a simple variant of it.

\subsubsection*{{Proof of Theorem \ref{thm_general_paraproduct}}} 
Let us first consider the case of strictly decreasing sequences $\lbrace a_j\rbrace_{j\in\mathbb{N}_0}, \, \lbrace b_j\rbrace_{j\in\mathbb{N}_0}$ satisfying hypothesis (Hyp 1).
For $j\in\mathbb{N}$ we define the frequency projections
\begin{align*}
    \widehat{\Delta^a_{j}f}=\mathbf{1}_{[a_{j+1},a_j)}\widehat{f}~, \quad  \widehat{\Delta^b_{j}g}=\mathbf{1}_{[b_j,b_0)}\widehat{g}~, \quad  \widehat{\Delta_{j}^c h}=\mathbf{1}_{-[a_{j+1},a_j)-[b_j,b_0)}\widehat{h}~. 
\end{align*}
With this notation, we rewrite the multiplier operator  $B_{m_{\mathtt{a,b}}}$ as 
$$B_{m_{\mathtt{a,b}}}(f,g)(x)= \sum_{j\in\mathbb{N}} \Delta^a_jf(x)\Delta^b_jg(x)~. $$
Let $h\in L^{p_3}(\mathbb{R})$ with $\Vert h \Vert_{p_3}=1$. Recall that, a bilinear multiplier acting on a pair of functions with frequency supports in the intervals $A$ and $B$, respectively, results in a function with frequency support in the interval $-A-B$. Hence we have that
$$\int_{\mathbb{R}} B_{m_{\mathtt{a,b}}}(f,g)(x) h(x) \dd x = \int_{\mathbb{R}} \sum_{j\in\mathbb{N}} \Delta^a_jf(x)\Delta^b_jg(x) \Delta^c_j h(x) \dd x~, $$
and by applying Hölder's inequality twice we obtain
\begin{align*}
\int_{\mathbb{R}} B_{m_{\mathtt{a,b}}}(f,g)(x) & h(x) \dd x  \\
& \leq \Vert \lbrace \Delta^a_j f \rbrace_{j\in\mathbb{N}} \Vert_{L^{p_1}(\ell^2)} \Vert \lbrace \Delta^b_j g \rbrace_{j\in\mathbb{N}} \Vert_{L^{p_2}(\ell^\infty)} \Vert \lbrace \Delta^c_j h \rbrace_{j\in\mathbb{N}} \Vert_{L^{p_3}(\ell^2)}~,
\end{align*}
where $\tfrac{1}{p_1}+\tfrac{1}{p_2}+\tfrac{1}{p_3}=1$. To bound the second term, observe that the pointwise bound $|\Delta^b_j g(x)|\leq 2 \mathscr{C}g(x)$ holds with $\mathscr{C}$ the Carleson--Hunt operator
$$\mathscr{C}g(x):= \sup_{N\in\mathbb{R}} \bigg\vert \int_{-\infty}^N \widehat{g}(\xi) e^{2\pi i \xi x}\dd \xi\bigg\vert~.$$
Hence, by Carleson--Hunt theorem, it follows that $\Vert \lbrace \Delta^b_j g \rbrace_{j\in\mathbb{N}} \Vert_{L^{p_2}(\ell^\infty)} \lesssim_{p_2} \Vert g \Vert_{p_2}$ for $1<p_2<\infty$.

To bound the first factor we can appeal to Rubio de Francia's square function, as the intervals $[a_{j+1}, a_j)$ are pairwise disjoint, obtaining that  $ \Vert \lbrace \Delta^a_j f \rbrace_{j\in\mathbb{N}}\Vert_{L^{p_1}(\ell^2)} \lesssim_{p_1} \Vert f \Vert_{p_1}$ for all $2\leq p_1 <\infty$. If the sequence $\lbrace a_j \rbrace_{j\in\mathbb{N}_0}$ is Littlewood--Paley then such bound can be extended to all $1<p_1 <\infty$ with the implicit constant now depending also on the sequence $\lbrace a_j \rbrace_{j\in\mathbb{N}_0}$, see e.g. \cite{Bou89,Ba21} for the behaviour of the constant for $p\in (1,2)$ in the lacunary and high-order lacunary cases, respectively. In view of the hypothesis on the collection of intervals $\mathcal{I}_{\mathtt{a,b}}$, the third term can be treated by similar arguments, obtaining that $ \Vert \lbrace \Delta^a_j f \rbrace_{j\in\mathbb{N}}\Vert_{L^{p_3}(\ell^2)} \lesssim_{p_3} n(\mathtt{a,b}) \Vert f \Vert_{p_3}$ for $2\leq p_3 <\infty$ which can be extended to all $1<p_3<\infty$ if the subcollections are Littlewood--Paley up to accounting for the dependence of the implicit constant on the sequences $\lbrace a_j \rbrace_{j\in\mathbb{N}_0}$, $\lbrace b_j \rbrace_{j\in\mathbb{N}_0}$.
By duality, the proof of Theorem \ref{thm_general_paraproduct} under hypothesis (Hyp 1) is completed.

The case of strictly decreasing sequences $\lbrace a_j\rbrace_{j\in\mathbb{N}_0}, \, \lbrace b_j\rbrace_{j\in\mathbb{N}_0}$ satisfying hypothesis (Hyp 2) can be addressed as follows.
Consider the multiplier symbol
$$\sum_{j\in\mathbb{N}_0} \mathbf{1}_{[a_{j+1},a_j)}(\xi)\mathbf{1}_{(b_\infty,b_j)}(\eta)~.$$
The corresponding multiplier operator can be seen to be bounded in the range of exponents claimed in the statement via arguments similar to those utilized above. Let $a_\infty=\lim_{j\rightarrow  \infty} a_j$, $a_\infty\in\mathbb{R}\cup\lbrace - \infty\rbrace$. By rewriting
$$m_{\mathtt{a,b}}(\xi,\eta)= \mathbf{1}_{\lbrace a_\infty < \xi < a_0\rbrace}(\xi)\mathbf{1}_{\lbrace b_\infty <\eta < b_0 \rbrace}(\eta)-\sum_{j\in\mathbb{N}_0} \mathbf{1}_{[a_{j+1},a_j)}(\xi)\mathbf{1}_{(b_\infty,b_j)}(\eta)~, $$
the claimed boundedness results for $B_{m_{\mathtt{a,b}}}$ under hypothesis (Hyp 2) follows. 
\qed

\vspace{0.2cm}

For later use, we record the following simple variant of Theorem \ref{thm_general_paraproduct}  which considers the case of strictly increasing sequences. Let ${{\mathtt{u}}}:=\lbrace {u}_j\rbrace_{j\in\mathbb{N}_0}, \; {\mathtt{v}}:=\lbrace {v}_j \rbrace_{j\in\mathbb{N}_0}$ be two sequences of strictly increasing real numbers. 
Our key assumption on these sequences reads as follows: the collection of intervals 
$\lbrace -({u}_0, {u}_{j}] - [{v}_{j}, {v}_{j+1}) \,\rbrace_{j\in\mathbb{N}}~$
can be split into a finite number ${n}({\mathtt{u}},{\mathtt{v}})$ of subcollections whose intervals are pairwise disjoint.
\begin{thm}\label{thm_general_paraproduct_bis}
    Let $\lbrace {u}_j\rbrace_{j\in\mathbb{N}_0}, \;\lbrace {v}_j \rbrace_{j\in\mathbb{N}_0}$ be as above. Consider the multiplier
    \begin{align*}
   \widetilde{m}_{\mathtt{{u},{v}}}(\xi,\eta) := \sum_{j\in\mathbb{N}} \mathbf{1}_{({u}_{0},{u}_j]}(\xi)\mathbf{1}_{[{v}_j,{v}_{j+1})}(\eta)
    \end{align*}
    and the associated bilinear multiplier operator 
    \begin{align*}
        B_{\widetilde{m}_{\mathtt{{u},{v}}}}(f,g)(x):=\int\int_{\mathbb{R}^2} \widetilde{m}_{\mathtt{{u},{v}}}(\xi,\eta)\widehat{f}(\xi)\widehat{g}(\eta)e^{2\pi i (\xi+\eta)x} \dd \xi \dd \eta~.
    \end{align*}
    The following estimate
    \begin{equation}\label{main_estimate_general_paraproduct_bis}
        \Vert B_{\widetilde{m}_{\mathtt{{u},{v}}}}(f,g) \Vert_{p_3'} \leq C(p_1,p_2,\mathtt{u,v}) \Vert f \Vert_{p_1} \Vert g \Vert_{p_2}~,
    \end{equation}
    holds for all Schwartz functions $f,g\in\mathscr{S}(\mathbb{R})$ whenever
    \begin{equation*}
       1<p_1<\infty ~, \quad 2\leq p_2, p_3< \infty~, \quad \text{and} \quad \frac{1}{p_1}+\frac{1}{p_2}+\frac{1}{p_3} =1~,
    \end{equation*}
    and, in such range of exponents, $C(p_1,p_2,\mathtt{u,v})=C(p_1,p_2,n(\mathtt{u,v}))$.
    Moreover, if the sequence $\lbrace {v}_j \rbrace_{j\in\mathbb{N}_0}$ is Littlewood--Paley then \eqref{main_estimate_general_paraproduct_bis} holds whenever
    \begin{equation*}
        1 < p_1, p_2 < \infty~, \quad 2\leq p_3 <\infty~, \quad \text{and} \quad \frac{1}{p_1}+\frac{1}{p_2}+\frac{1}{p_3} =1~.
    \end{equation*}
    If, in addition, the collection of intervals $\lbrace -({u}_0, {u}_{j}] - [{v}_{j}, {v}_{j+1}) \,\rbrace_{j\in\mathbb{N}}~$ can be split into a finite number of subcollections which are also Littlewood--Paley then \eqref{main_estimate_general_paraproduct_bis} holds whenever 
    \begin{equation*}
        1 < p_1, p_2, p_3 < \infty~,\quad \text{and} \quad \frac{1}{p_1}+\frac{1}{p_2}+\frac{1}{p_3} =1~.
    \end{equation*}
\end{thm}

\section{Proof of Theorem \ref{thm_1_over_x_power} and of Theorem \ref{thm_power_positive}}\label{section_proof_1_over_x_power}

\subsubsection*{Proof of Theorem \ref{thm_1_over_x_power}}  Let $m_{\Gamma_0}$
be the characteristic function of the set
$$\Gamma_0:=\lbrace (\xi,\eta): \, \xi < - c^{1/(c+1)}, \,  \gamma(\xi)\leq \eta < c^{-c/(c+1)} \rbrace~.$$
Observe that $\gamma'(- c^{1/(c+1)})=1$ and $\gamma(-c^{1/(c+1)})=c^{-c/(c+1)}$. In particular, the cut-off at $-c^{1/(c+1)}$ ensures that $0<\gamma'<1$ on $(-\infty, - c^{1/(c+1)})$.
To prove Theorem \ref{thm_1_over_x_power} it is enough to show that the bilinear multiplier $B_{m_{\Gamma_0}}$ with symbol  $m_{\Gamma_0}$ and defined as in \eqref{multiplier_operator} is a bounded $L^{p_1}\times L^{p_2}\rightarrow L^{p_3'}$ operator for all $p_1,\, p_2,\, p_3$ in the local $L^2$ range \eqref{local_L2}. In fact, the set $\lbrace (\xi,\eta): \,  \xi < 0, \, \eta \geq \gamma(\xi) \rbrace$ can be decomposed as the union
\begin{align*}
    \; \Gamma_0  & \cup \lbrace (\xi, \eta):  \xi \leq  - c^{1/(c+1)},\, c^{-c/(c+1)} \leq \eta  \rbrace\\
   & \; \cup \lbrace (\xi,\eta): \, c^{-c/(c+1)} < \eta , \,  - c^{1/(c+1)} \leq \xi < - \eta^{-1/c} \rbrace ~.
\end{align*}
The boundedness of the bilinear multiplier operator associated with the second set can be seen to follow, for example, from the $L^p$ boundedness properties of the (modulated) Hilbert transform, while the bilinear multiplier operator associated with the third set can be treated by arguments similar to those that we are about to illustrate here by inverting the role of $\xi$ and $\eta$.

We further decompose the multiplier corresponding to $\Gamma_0$ into the sum of a staircase paraproduct and some boundary terms.
To achieve this, we start by constructing a sequence $\lbrace a_j \rbrace_{j\in\mathbb{N}_0}$ as follows. 
For each $j\in\mathbb{N}_0$ we choose $a_j$ so that 
$2^{-j}=\gamma'(a_j)$, i.e. $a_j=-(c2^{j})^{1/(c+1)}$. In particular, $a_0=-c^{1/(c+1)}$ and the sequence is strictly decreasing.
One can easily check that for this choice of the sequence $\lbrace a_j\rbrace_{j\in\mathbb{N}_0}$ we have that over each interval $[a_{j+1},a_j)$ the slope of the curve $\gamma$ is bounded between $\alpha_j$ and $2\alpha_j$ for some $0<\alpha_j<1$.
Moreover, we observe that the just constructed sequence $\lbrace a_j \rbrace_{j\in\mathbb{N}_0}$ is lacunary as $|a_{j+1}|/|a_j|=2^{1/(c+1)}$ for all $j\in\mathbb{N}_0$ and, therefore, it is a Littlewood--Paley sequence.
Then, we define $\lbrace b_j \rbrace_{j\in\mathbb{N}_0}$ to be the sequence of elements of the form $b_j:=\gamma(a_j)$. The sequence $\lbrace b_j \rbrace_{j\in\mathbb{N}_0}$ is therefore strictly decreasing and for all $j\in \mathbb{N}_0$ we have $0<b_j\leq c^{-c/(c+1)}$.
Let $\mathcal{I}_{\mathtt{a,b}}$ be the collection of intervals $\lbrace -[a_{j+1},a_j)-[b_j,b_0) :\, j\in\mathbb{N}\, \rbrace$. Note that for each $j\in\mathbb{N}$ we have
 $$-[a_{j+1},a_j)-[b_j,b_0) \subseteq (|a_j|-c^{-c/(c+1)}, |a_{j+1}|)~.$$
 Lacunarity of the $a_j$'s guarantees that there exists a ${j}_0$ such that $|a_j-a_{j-1}|> c^{-c/(c+1)}$ for all $j\geq j_0$.
 Hence, the collection of intervals  
 $\mathcal{I}_{\mathtt{a,b}}$ can be split into a finite number of subcollections so that the intervals in each subcollection are pairwise disjoint and, in view of the lacunarity of the sequence $\lbrace a_j \rbrace_{j\in\mathbb{N}_0}$, such that each subcollection is Littlewood--Paley.
 It follows from Theorem \ref{thm_general_paraproduct} that the bilinear multiplier operator $B_{m_{par}}$ with symbol
$$m_{\mathtt{a,b}}(\xi,\eta):= \sum_{j\in\mathbb{N}} \mathbf{1}_{[a_{j+1}, a_j)}(\xi)\mathbf{1}_{[b_j,b_0)}(\eta)$$
is a bounded $L^{p_1}\times L^{p_2} \rightarrow L^{p_3'}$ operator for all $$1<\, p_1,\, p_2, \, p_3\,  <\infty \quad \text{such that} \quad \frac{1}{p_1}+\frac{1}{p_2}+\frac{1}{p_3}=1~.$$

For each $j\in\mathbb{N}_0$ define
$$m_j(\xi, \eta):= \mathbf{1}_{[a_{j+1},a_j)}(\xi)\mathbf{1}_{[\gamma(\xi), b_j)}(\eta)~.$$
Then we can decompose $m=\mathbf{1}_{\Gamma_0}$ as 
$$m(\xi, \eta)= m_{\mathtt{a,b}}(\xi,\eta)+ \sum_{j\in\mathbb{N}_0} m_j(\xi,\eta)~.$$
By construction, the slope of the curved boundary line of each of the sets
$$\lbrace (\xi,\eta) : a_{j+1}\leq \xi < a_j, \, \gamma(\xi) \leq \eta < b_j \rbrace$$
is between $\alpha_j$ and $2\alpha_j$ for some $0<\alpha_j<1$. 
Because of this, each bilinear multiplier operator $B_{m_j}$ can be seen to be individually and -- crucially-- uniformly bounded in the local $L^2$ range \eqref{local_L2} relying, for example, on the analysis in \cite[Section~6]{GL06} combined with previous results from \cite{GL04,L06}. See also \cite{M00}, and see also \cite{ST23} where this same crucial property has been used.

We conclude the proof of the theorem by arguing as in \cite[Lemma~1]{GL06}, see also \cite[Proof~of~Corollary~1.2]{ST23} for a similar argument.
Define the frequency projections
\begin{align*}
    \widehat{\Delta^a_{j}f}:=\mathbf{1}_{[a_{j+1},a_j)}\widehat{f}~, \quad  \widehat{\Delta^b_{j}g}:=\mathbf{1}_{[b_{j+1},b_j)}\widehat{g}~, \quad  \widehat{\Delta_{j}^c h}:=\mathbf{1}_{-[a_{j+1},a_j)-[b_{j+1},b_j)}\widehat{h}~,
\end{align*}
and observe that we can rewrite
$$\sum_{j\in\mathbb{N}_0}B_{m_j}(f,g)(x)= \sum_{j\in\mathbb{N}_0}\Delta_{j}^c B_{m_j}(\Delta^a_{j} f, \Delta^b_{j}g)(x)~.$$
It follows from the discussion above that, up to splitting into a finite number of subcollections, the intervals in the collection $\lbrace -[a_{j+1},a_j)-[b_{j+1},b_j) \rbrace_{j\in\mathbb{N}_0}$ are pairwise disjoint.

Let $h\in L^{p_3}(\mathbb{R})$ with $\Vert h \Vert_{L^{p_3}}=1$. Using the properties of the frequency supports, Hölder's inequality (twice), Rubio de Francia's square function, and inclusion properties of  $\ell^p$ spaces, we obtain
\begin{align*}
    \int_{\mathbb{R}} \sum_{j\in\mathbb{N}_0}\Delta_{j}^c B_{m_j}(\Delta^a_{j} f, \Delta^b_{j}g)(x) & h(x) \dd x  = \int_{\mathbb{R}} \sum_{j\in\mathbb{N}_0} B_{m_j}(\Delta^a_{j} f, \Delta^b_{j}g)(x) \Delta_{j}^c  h(x) \dd x \\
    & \leq \Vert \lbrace B_{m_j}(\Delta_j^a f , \Delta^b_j g)\rbrace_{j\in\mathbb{N}_0} \Vert_{L^{p_3'}(\ell^2)} \Vert \lbrace \Delta_j^ch \rbrace_{j\in\mathbb{N}_0} \Vert_{L^{p_3}(\ell^2)} \\
   &  \lesssim \Vert \lbrace B_{m_j}(\Delta_j^a f , \Delta^b_j g)\rbrace_{j\in\mathbb{N}_0} \Vert_{L^{p_3'}(\ell^{2})} \Vert h \Vert_{L^{p_3}} \\
   & \lesssim \Vert \lbrace B_{m_j}(\Delta_j^a f , \Delta^b_j g)\rbrace_{j\in\mathbb{N}_0} \Vert_{L^{p_3'}(\ell^{p_3'})}~.
\end{align*}
Using the aforementioned uniform boundedness of the $B_{m_j}$'s, Hölder's inequality, the inclusion properties of $\ell^p$ spaces, and Rubio de Francia's square function we see that
\begin{align*}
    \Vert \lbrace B_{m_j}(\Delta_j^a f , \Delta^b_j g)\rbrace_{j\in\mathbb{N}_0} \Vert_{L^{p_3'}(\ell^{p_3'})} & = \bigg( \sum_{j\in\mathbb{N}_0} \Vert B_{m_j}(\Delta_j^a f , \Delta^b_j g) \Vert_{p_3'}^{p_3'}\bigg)^{1/p_3'} \\
    & \lesssim  \bigg( \sum_{j\in\mathbb{N}_0} \Vert \Delta_j^a f \Vert_{p_1}^{p_3'} \Vert \Delta_j^b g\Vert_{p_2}^{p_3'}\bigg)^{1/p_3'} \\
    &  \lesssim \bigg( \sum_{j\in\mathbb{N}_0} \Vert \Delta_j^a f \Vert_{p_1}^{p_1} \bigg)^{1/p_1} \bigg( \sum_{j\in\mathbb{N}_0}\Vert \Delta_j^b g\Vert_{p_2}^{p_2}\bigg)^{1/p_2} \\
    & \lesssim \Vert \lbrace \Delta_j^a f\rbrace_{j\in\mathbb{N}_0} \Vert_{L^{p_1}(\ell^{p_1})} \Vert \lbrace \Delta_j^b g\rbrace_{j\in\mathbb{N}_0} \Vert_{L^{p_2}(\ell^{p_2})} \\
    & \lesssim \Vert \lbrace \Delta_j^a f\rbrace_{j\in\mathbb{N}_0} \Vert_{L^{p_1}(\ell^2)} \Vert \lbrace \Delta_j^b g\rbrace_{j\in\mathbb{N}_0} \Vert_{L^{p_2}(\ell^2)}\\
    & \lesssim \Vert f \Vert_{L^{p_1}}\Vert h \Vert_{L^{p_2}}~,
\end{align*}
hence concluding the proof of Theorem \ref{thm_1_over_x_power}.
\qed

\subsubsection*{Proof of Theorem \ref{thm_power_positive}} Let $m_{\Gamma_0}$
be the characteristic function of the set
$$\Gamma_0:=\lbrace (\xi,\eta):\, \xi\in(0,\tfrac{1}{\sqrt{3}}), \, \gamma(\xi) \leq \eta < \gamma(\tfrac{1}{\sqrt{3}}) \rbrace~.$$
To prove Theorem \ref{thm_power_positive} it is enough to show that the bilinear multiplier $B_{m_{\Gamma_0}}$ with symbol  $m_{\Gamma_0}$ and defined as in \eqref{multiplier_operator} is a bounded $L^{p_1}\times L^{p_2}\rightarrow L^{p_3'}$ operator for all $p_1,\, p_2,\, p_3$ in the local $L^2$ range \eqref{local_L2}. We further decompose $B_{m_{\Gamma_0}}$ into the sum of a staircase paraproduct and some boundary terms. To this end, let  $\lbrace a_j\rbrace_{j\in\mathbb{N}}$ be the sequence of positive real numbers satisfying the condition $\gamma'(a_j)=2^{-j}$ for all $j\in\mathbb{N}$. In particular, we have that $a_j=\tfrac{1}{\sqrt{2^{2j}-1}}$ and the sequence is decreasing.  Let $\lbrace b_j\rbrace_{j\in\mathbb{N}}$ be the sequence of elements of the form $b_j=\gamma(a_j)$. The sequence is strictly decreasing and $\lim_{j\rightarrow\infty}b_j=1$. Let $\mathcal{J}_{\mathtt{a,b}}$ be the collection of intervals $\lbrace -[a_{j+1},a_j)-(b_\infty,b_j) \rbrace_{j\in\mathbb{N}}$. Using the crude estimate $b_j<a_j+1$ we see that for each $j\in\mathbb{N}$
$$ -[a_{j+1},a_j)-(b_\infty,b_j) \subseteq (-2a_j -1, -a_{j+1}-1]~.$$
Hence the collection of intervals $\mathcal{J}_{\mathtt{a,b}}$ can be split into two subcollections whose intervals are pairwise disjoint. Boundedness of the paraproduct operator associated with the sequences  $\lbrace a_j\rbrace_{j\in\mathbb{N}}$, $\lbrace b_j\rbrace_{j\in\mathbb{N}_0}$ follows from Theorem \ref{thm_general_paraproduct}. The choice of the sequence $\lbrace a_j \rbrace_{j\in\mathbb{N}_0}$ guarantees that the passage to bounds for the multiplier in the local $L^2$ range of exponents \eqref{local_L2} can be achieved by the same arguments used in the proof of Theorem \ref{thm_1_over_x_power}.
\qed

\section{Further examples and generalizations}\label{section_examples}

\subsection*{Epigraph of exponential curves}\label{subsec:exp}
In this subsection, we briefly discuss the boundedness of a paraproduct associated with the exponential curve $\eta = 2^\xi$. 
Let $m_{\mathtt{par}}$ be the multiplier symbol
\begin{align*}
    m_{\mathtt{par}}(\xi,\eta)  := & m_1(\xi,\eta)+m_2(\xi,\eta)+m_3(\xi,\eta) \\
     := & \sum_{j\in\mathbb{N}_0} \mathbf{1}_{[-(j+1),-j)}(\xi)\mathbf{1}_{[2^{-j},1)}(\eta) + \sum_{j\in\mathbb{N}} \mathbf{1}_{(0,j)}(\xi)\mathbf{1}_{[2^j, 2^{j+1})}(\eta) \\
     & + \mathbf{1}_{\lbrace \xi \leq 0\rbrace}(\xi)\mathbf{1}_{\lbrace \eta\geq 1 \rbrace}(\eta)~.
\end{align*}
Hence, $m_{\mathtt{par}}$ is the indicator function of a staircase set inscribed by the exponential curve $\eta=2^\xi$. Let $B_{m_\mathtt{par}}$, $B_{m_1}$, $B_{m_2}$, $B_{m_3}$ be the bilinear multiplier operators associated with $m_{\mathtt{par}}, \, m_1,\,m_2,\, m_3$, respectively, defined as in \eqref{multiplier_operator}. $B_{m_1}$ is a paraproduct covered by Theorem \ref{thm_general_paraproduct}. In particular, $B_{m_1}$ is a bounded $L^{p_1}\times L^{p_2}\rightarrow L^{p_3'}$ operator for all
$$2\leq \, p_1, \, p_3 \, <\infty \,, \quad 1<p_2<\infty \,, \qquad \frac{1}{p_1}+\frac{1}{p_2}+\frac{1}{p_3}=1~.$$
Note that boundedness of $B_{m_1}$ in the open local $L^2$ region
$2 < \, p_1, \,p_2, \, p_3 \, <\infty$, $\tfrac{1}{p_1}+\tfrac{1}{p_2}+\tfrac{1}{p_3}=1~,$ was established first in \cite{ST23} via different arguments. The bilinear operator $B_{m_2}$ is a paraproduct covered by Theorem \ref{thm_general_paraproduct_bis} and it is a bounded 
$L^{p_1}\times L^{p_2}\rightarrow L^{p_3'}$ operator for all
$$1 < \, p_1, \, p_2 \, <\infty \,, \quad 2\leq p_3<\infty \,, \qquad \frac{1}{p_1}+\frac{1}{p_2}+\frac{1}{p_3}=1~.$$ Boundedness of the bilinear multiplier operator $B_{m_3}$ for all $1<\, p_1,\, p_2\, p_3 <\infty$, $\tfrac{1}{p_1}+\tfrac{1}{p_2}+\tfrac{1}{p_3}=1~,$ is standard. As a consequence, the paraproduct $B_{m_{\mathtt{par}}}$ is a bounded $L^{p_1}\times L^{p_2}\rightarrow L^{p_3'}$ operator for all
$$2\leq \, p_1, \, p_3 \, <\infty \,, \quad 1<p_2<\infty \,, \qquad \frac{1}{p_1}+\frac{1}{p_2}+\frac{1}{p_3}=1~.$$

In the local $L^2$ range of exponents \eqref{local_L2}, the passage from estimates for the paraproduct $B_{m_\mathtt{par}}$ to bounds for the multiplier $B_m$, defined as in \eqref{multiplier_operator} with $m$ the characteristic function of the set
$$\lbrace (\xi,\eta): \, \xi\in\mathbb{R}, \, \eta \geq 2^\xi \rbrace~,$$
can be achieved by the procedure outlined in \cite{ST23} relying on the results from \cite{GL06,M00} and observing that $B_{m}((\mathbf{1}_{(0,\infty)}\widehat{f})^\vee,g)-B_{m_2}(f,g)$ can be dealt with as in \cite{ST23} by reversing the role of $\xi$ and $\eta$.

\subsection*{Epigraph of polygonal curves} Let $\lbrace a_j \rbrace_{j\in\mathbb{N}_0}$, $\lbrace b_j \rbrace_{j\in\mathbb{N}_0}$ be two sequences of strictly decreasing real numbers satisfying hypothesis (Hyp 2). Let $\gamma$ be the piecewise linear curve of vertexes $\lbrace (a_j,b_j) \rbrace_{j\in\mathbb{N}_0}$ and define $m$ to be the indicator function of the epigraph of $\gamma$. Note in passing that the polygonal curve $\gamma$ may not be convex. The corresponding bilinear multiplier operator $B_m$ defined as in \eqref{multiplier_operator} is a bounded $L^{p_1}\times L^{p_2}\rightarrow L^{p_3'}$ operator for all $p_1,\, p_2,\, p_3$ in the local $L^2$ range \eqref{local_L2}.  In fact, boundedness for the staircase paraproduct associated with the sequences $\lbrace a_j \rbrace_{j\in\mathbb{N}_0}$, $\lbrace b_j \rbrace_{j\in\mathbb{N}_0}$ is established by Theorem \ref{thm_general_paraproduct}. Then, the passage to bounds for the multiplier in the local $L^2$ range of exponents \eqref{local_L2} can be achieved as follows. Let $H_{\ell_j}$ be the bilinear multiplier operator whose symbol is the indicator function of the triangle of vertexes $(a_j,b_j)$, $(a_{j+1},b_j)$, and $(a_{j+1},b_{j+1})$.
Let $\Delta_j^a f:= \mathbf{1}_{[a_{j+1},a_j)}\widehat{f}$, $\Delta_j^b g:= \mathbf{1}_{[b_{j+1},b_j)}\widehat{g} $.  Then, in view of the frequency support properties, uniform bounds for the bilinear Hilbert transform from \cite{GL04,L06,GL06}, and Rubio de Francia's square function, arguing as in Section \ref{section_proof_1_over_x_power} one can easily check that for all triples of exponents $p_1,\, p_2,\, p_3$ in the local $L^2$ range \eqref{local_L2} we have
\begin{align*}
    \bigg\Vert \sum_{j\in\mathbb{N}_0} H_{\ell_j} (\Delta_j^af, \Delta_j^b g) \bigg\Vert_{L^{p_3'}(\mathbb{R})} & \lesssim \Vert \lbrace H_{\ell_j}(\Delta_j^af, \Delta_j^b g) \rbrace_{j\in\mathbb{N}_0} \Vert_{L^{p_3'}(\ell^2)} \\
    & \lesssim \Vert \lbrace \Delta_j^af \rbrace_{j\in\mathbb{N}_0} \Vert_{L^{p_1}(\ell^2)} \Vert \lbrace \Delta_j^b g \rbrace_{j\in\mathbb{N}_0} \Vert_{L^{p_2}(\ell^2)}\\
    & \lesssim \Vert f \Vert_{L^{p_1}(\mathbb{R})} \Vert g \Vert_{L^{p_2}(\mathbb{R})}~.
\end{align*}

\subsection*{Convex sequences} Let $\gamma\in\mathscr{C}^1$ be an increasing strictly convex function on $(-\infty,c)$ such that $0 < \gamma'<1$ on $(-\infty,a_0)$, $-\infty<a_0<c$. For simplicity, we may assume that $\gamma'(a_0)=1$, $a_0=0$, and $\gamma(a_0)=0$. Let $\lbrace a_j\rbrace_{j\in\mathbb{N}_0}$ be the sequence of elements satisfying the condition $\gamma'(a_j)=2^{-j}$. In particular, $\gamma'(a_0)=1$ and the sequence is decreasing. Recall that a sequence of strictly increasing, positive real numbers $\lbrace \alpha_j\rbrace_{j\in\mathbb{N}_0}$ is said to be convex if $\alpha_j \leq \tfrac{1}{2}(\alpha_{j-1}+\alpha_{j+1})$, for all $j\in\mathbb{N}$. This is equivalent to the condition
$$\frac{\alpha_{j+1}-\alpha_j}{\alpha_j-\alpha_{j-1}}\geq 1 \qquad \text{for all}\; j\in\mathbb{N}~.$$
We say that our sequence $\lbrace a_j\rbrace_{j\in\mathbb{N}_0}$ is convex if $\lbrace |a_j|\rbrace_{j\in\mathbb{N}_0}$ is convex. Let $\lbrace b_j\rbrace_{j\in\mathbb{N}_0}$ be the sequence of elements of the form $b_j=\gamma(a_j)$. The convexity of the curve $\gamma$, together with the definition of the $a_j's$ and $b_j$'s, implies that 
$$\frac{|b_j-b_{j+1}|}{|a_j-a_{j+1}|}\leq 2^{-j}~.$$
With this observation and by convexity of the sequence $\lbrace a_j\rbrace_{j\in\mathbb{N}_0}$, we obtain that
$$|b_j|=\sum_{k=1}^{j}|b_{k}-b_{k-1}|\leq \sum_{k=1}^{j}2^{-(k-1)}|a_{k}-a_{k-1}|\leq 2|a_{j}-a_{j-1}|~.$$
As a consequence, we have that for all $j\in\mathbb{N}$
\begin{align*}
    -(a_{j+1},a_j]-(b_j,b_0] & = [|a_j|, |a_{j+1}|) +[0,|b_{j}|)\\
    & \subset [|a_j|, |a_{j+1}|+2|a_j-a_{j-1}|)~.
\end{align*}
Hence, the collection of intervals $\lbrace -[a_{j+1},a_j) -[b_{j},b_0) \, \rbrace_{j\in\mathbb{N}}$
can be split into three subcollections whose intervals are pairwise disjoint. It follows from Theorem \ref{thm_general_paraproduct} that for such sequences $\mathtt{a}:=\lbrace a_j\rbrace_{j\in\mathbb{N}_0}$ and $\mathtt{b}:=\lbrace b_j\rbrace_{j\in\mathbb{N}_0}$ the associated paraproduct operator $B_{m_{\mathtt{a,b}}}$ defined as in \eqref{defi_B_m_a_b} is a bounded $L^{p_1}\times L^{p_2}\rightarrow L^{p_3'}$ operator for all exponents 
$$ 2\leq p_1, p_3 < \infty~, \quad 1 < p_2 <\infty~, \quad \text{and} \quad \frac{1}{p_1}+\frac{1}{p_2}+\frac{1}{p_3} =1~.$$
Let $m$ be the indicator function of the epigraph of $\gamma$ restricted to $(-\infty,a_0)$ and let $B_m$ be the associated multiplier operator defined as in \eqref{multiplier_operator}.
The choice of the sequence $\lbrace a_j \rbrace_{j\in\mathbb{N}_0}$ and its convexity guarantee that the passage to bounds for the multiplier in the local $L^2$ range of exponents \eqref{local_L2} can be achieved by the same arguments detailed in Section \ref{section_proof_1_over_x_power}. A simple example of a function $\gamma$ for which the sequence $\lbrace (\gamma')^{-1}(2^{-j}) \rbrace_{j\in\mathbb{N}_0}$ is convex is the exponential function $\gamma(\xi)=2^\xi$.

\subsection*{Lacunary sequences}
Let $\gamma\in\mathscr{C}^1$ be an increasing strictly convex function on $(-\infty,c)$ such that $0 < \gamma'<1$ on $(-\infty,a_0)$,  $-\infty< a_0<c$. For simplicity, we may assume that $\gamma'(a_0)=1$, $a_0=0$, and $\gamma(a_0)=0$. We define $\lbrace a_j\rbrace_{j\in\mathbb{N}_0}$ to be the sequence of elements satisfying the condition $\gamma'(a_j)=2^{-j}$. In particular, we have that $\gamma'(a_0)=1$ and the sequence is decreasing. Recall that a sequence of positive real numbers $\lbrace \alpha_j\rbrace_{j\in\mathbb{N}_0}$ is said to be lacunary if there exists $q>1$ such that $\alpha_{j+1}/\alpha_j\geq q$ for all $j\in\mathbb{N}_0$. We say that our sequence $\lbrace a_j\rbrace_{j\in\mathbb{N}_0}$ is lacunary if $\lbrace |a_j|\rbrace_{j\in\mathbb{N}_0}$ is lacunary. Let $\lbrace b_j\rbrace_{j\in\mathbb{N}_0}$ be the sequence of elements of the form $b_j=\gamma(a_j)$.  Using the crude estimate $|b_j|<|a_j|$ we have that for all $j\in\mathbb{N}$
\begin{align*}
    -(a_{j+1},a_j]-(b_j,b_0] & = [|a_j|, |a_{j+1}|) +[0,|b_{j}|)\\
    & \subset [|a_j|, |a_{j+1}|+|a_j|)~.
\end{align*}
Hence, the collection of intervals $\lbrace -[a_{j+1},a_j) -[-b_{j},b_0) \, \rbrace_{j\in\mathbb{N}}$
can be split into 
a finite number $n(\mathtt{a,b})=n(q)$ of subcollections whose intervals are pairwise disjoint and, thanks to the lacunarity of the $a_j$'s, such that these subcollections are Littlewood--Paley. It follows from Theorem \ref{thm_general_paraproduct} that for such sequences $\mathtt{a}:=\lbrace a_j\rbrace_{j\in\mathbb{N}_0}$ and $\mathtt{b}:=\lbrace b_j\rbrace_{j\in\mathbb{N}_0}$ the associated paraproduct operator $B_{m_{\mathtt{a,b}}}$ defined as in \eqref{defi_B_m_a_b} is a bounded $L^{p_1}\times L^{p_2}\rightarrow L^{p_3'}$ operator for all exponents 
$$ 1 < p_1, p_2, p_3 < \infty~,  \quad \text{and} \quad \frac{1}{p_1}+\frac{1}{p_2}+\frac{1}{p_3} =1~.$$
Let $m$ be the indicator function of the epigraph of $\gamma$ restricted to $(-\infty,a_0)$ and let $B_m$ be the associated multiplier operator defined as in \eqref{multiplier_operator}.
The choice of the sequence $\lbrace a_j \rbrace_{j\in\mathbb{N}_0}$ and its lacunarity guarantee that the passage to bounds for the multiplier in the local $L^2$ range of exponents \eqref{local_L2} can be achieved by the same arguments outlined in Section \ref{section_proof_1_over_x_power}. 
 Simple examples of convex functions $\gamma$ for which the sequence $\lbrace (\gamma')^{-1}(2^{-j}) \rbrace_{j\in\mathbb{N}_0}$ is lacunary include the family of functions considered in Theorem \ref{thm_1_over_x_power}, rational functions like $\gamma(\xi)=\tfrac{\xi}{\xi+c}$, $c>0$, on $(-\infty,-c)$, the function $\gamma(\xi)=\arctan(\xi)$ on $(-\infty,0]$.

 \subsection*{Concave sequences} Let $\gamma\in\mathscr{C}^1$ be an increasing, strictly convex function on $(a_\infty,c)$ such that $0<\gamma'<1$ on $(a_\infty,a_0)$, $a_\infty<a_0<c $, and such that $\lim_{\xi\rightarrow a_\infty}\gamma(\xi)= b_\infty$ for some $b_\infty\in\mathbb{R}$. For simplicity, we may assume that $a_\infty=0$, $b_\infty=0$, and $\gamma'(a_0)=1$. Let $\lbrace a_j \rbrace_{j\in\mathbb{N}_0}$ be the sequence of elements satisfying the condition $\gamma'(a_j)=2^{-j}$. In particular, $\gamma'(a_0)=1$ and the sequence is strictly decreasing. Recall that a sequence of strictly increasing, positive real numbers is said to be concave if $\alpha_j\geq \tfrac{1}{2}(\alpha_{j+1}+\alpha_{j-1})$ for all $j\in\mathbb{N}$. This is equivalent to the condition 
 $$\frac{\alpha_{j}-\alpha_{j-1}}{\alpha_{j+1}-\alpha_{j}}\geq 1 \qquad \text{for all}\; j\in\mathbb{N}~.$$
 We say that our sequence $\lbrace a_j \rbrace_{j\in\mathbb{N}_0}$ is concave if $|a_j-a_{j-1}|\geq |a_{j+1}-a_j|$ for all $j\in\mathbb{N}$. Let $\lbrace b_j \rbrace_{j\in\mathbb{N}_0}$ be the sequence of elements of the form $b_j=\gamma(a_j)$. It follows from our hypothesis that this is a strictly decreasing sequence. Convexity of the curve $\gamma$ together with concavity of the sequence $\lbrace a_j \rbrace_{j\in\mathbb{N}_0}$ implies that
 $$|b_j|\leq \sum_{k=j}^\infty 2^{-k}|a_{k+1}-a_k|\leq 2 |a_j-a_{j+1}|~.$$
 It follows that for all $j\in\mathbb{N}_0$
 $$-[a_{j+1},a_j)-(b_\infty,b_j)=(-a_j-b_j,-a_{j+1}]\subseteq (-a_j-2|a_j-a_{j+1}|, -a_{j+1}]~.$$
 Therefore, the collection of intervals $\lbrace -[a_{j+1},a_j)-(b_\infty,b_j) \rbrace_{j\in\mathbb{N}_0}$ can be split into three subcollections whose intervals are pairwise disjoint. It follows from Theorem \ref{thm_general_paraproduct} that for such sequences $\mathtt{a}:=\lbrace a_j\rbrace_{j\in\mathbb{N}_0}$ and $\mathtt{b}:=\lbrace b_j\rbrace_{j\in\mathbb{N}_0}$ the associated paraproduct operator $B_{m_{\mathtt{a,b}}}$ defined as in \eqref{defi_B_m_a_b} is a bounded $L^{p_1}\times L^{p_2}\rightarrow L^{p_3'}$ operator for all exponents 
$$ 2\leq p_1, p_3 < \infty~, \quad 1 < p_2 <\infty~, \quad \text{and} \quad \frac{1}{p_1}+\frac{1}{p_2}+\frac{1}{p_3} =1~.$$
Let $m$ be the indicator function of the epigraph of $\gamma$ restricted to $(a_\infty,a_0)$ and let $B_m$ be the associated multiplier operator defined as in \eqref{multiplier_operator}.
The choice of the sequence $\lbrace a_j \rbrace_{j\in\mathbb{N}_0}$ guarantees that the passage to bounds for the multiplier in the local $L^2$ range of exponents \eqref{local_L2} can be achieved by the same arguments detailed in Section \ref{section_proof_1_over_x_power}. Simple examples of convex, increasing functions for which the sequence $\lbrace (\gamma')^{-1}(2^{-j}) \rbrace_{j\in\mathbb{N}_0}$ is concave are $\gamma(\xi)=\xi^c/c$, $c>0$, on $(0,1)$, the function considered in Theorem \ref{thm_power_positive}, $\gamma(\xi)=\sqrt{1+\xi^2}$ on $(0,\tfrac{1}{\sqrt{3}})$, and (up to a finite number of exceptions) $\gamma(\xi)= - \sqrt{1-\xi^2}$ on $(0,\tfrac{1}{\sqrt{2}})$.

\subsection*{Sufficient conditions}  We summarize some of the above in the following statement.
\begin{thm}
Let $\gamma$ be a strictly convex, increasing function satisfying one of the following:
   \begin{enumerate}[label=(\roman*)]
        \item[(i)] $\gamma$ is the restriction to $(-\infty,0)$ of a $\mathscr{C}^1$ function such that $0<\gamma'<1$ on $(-\infty,0)$, $\gamma'(0)=1$, and  the sequence $\lbrace (\gamma')^{-1} (2^{-j}) \rbrace_{j\in\mathbb{N}_0}$ is convex;
        \item[(ii)] $\gamma$ is the restriction to $(-\infty,0)$ of a $\mathscr{C}^1$ function such that $0<\gamma'<1$ on $(-\infty,0)$, $\gamma'(0)=1$, and  the sequence $\lbrace (\gamma')^{-1} (2^{-j}) \rbrace_{j\in\mathbb{N}_0}$ is lacunary;
        \item[(iii)] $\gamma$ is the restriction to $(0,1)$ of a $\mathscr{C}^1$ function such that $0<\gamma'<1$ on $(0,1)$, $\gamma(0)=0$, $\gamma'(1)=1$, and the sequence $\lbrace (\gamma')^{-1} (2^{-j}) \rbrace_{j\in\mathbb{N}_0}$ is concave.
    \end{enumerate}
    Let $m$ be the characteristic function of the epigraph of $\gamma$. Then the bilinear multiplier operator $B_m$ defined as in \eqref{multiplier_operator} satisfies the a priori bound \eqref{wanted_bound} for all triple of exponents $p_1,p_2,p_3$ in the local $L^2$ range \eqref{local_L2}.
    
\end{thm}

\section{Proof of Theorem \ref{thm_polygonal}}\label{section_polygonal}

Boundedness of the bilinear multiplier operator $B_m$ in the local $L^2$ range \eqref{local_L2} follows from the second example in Section \ref{section_examples}. Therefore, in this Section we focus on studying the boundedness of $B_m$ in the range of exponents
\begin{align}\label{range_beyond_L2}
    2< \, p_1,\, p_3 \, < \infty\, , \quad 1 < p_2 < 2\, , \quad \frac{1}{p_1}+\frac{1}{p_2}+\frac{1}{p_3}=1~.
\end{align}
To achieve the passage from estimates for the staircase paraproduct to bounds for our multiplier in the exponent range \eqref{range_beyond_L2} we rely on the approach developed by Demeter and Gautam \cite{DG12}. In particular, via a classical time-frequency discretization procedure one can reduce boundedness of the multiplier operator to boundedness of a discretized model sum. This can be treated by invoking \cite[Theorem~2.3]{DG12}, also inspired by the analysis in \cite{MTT02}. Here we detail the derivation of the discretized model sum for our multiplier. We follow very closely the steps and the notation of \cite{DG12} accounting for the minor modifications required by our setting.

Let $T_j, \, j\in\mathbb{N}_0$ be the triangle with vertices in $(a_j,b_j),$ $(a_{j+1},b_j),$ and $(a_{j+1},b_{j+1})$. We want to cover each triangle $T_j$ with a collection of Whitney rectangles. To achieve this, we introduce the collection ${\mathbf{S}}$ of squares  $S\subset \mathbb{R}^2$ with centers in the lattice $2^{j-10}\mathbb{Z}^2$, of side length $2^j$ for some $j\in\mathbb{Z}$, and satisfying the Whitney condition
\begin{align*}
    C_0\, S\cap\lbrace (\xi,\xi):\, \xi\in\mathbb{R} \rbrace & = \emptyset \\
    4 C_0\, S\cap\lbrace (\xi,\xi):\, \xi\in\mathbb{R} \rbrace & \neq \emptyset
\end{align*}
for a suitable, sufficiently large, constant $C_0$. Consider the family $\lbrace L_j^2 \rbrace_{j\in\mathbb{N}_0}$ of linear transformations on $\mathbb{R}^2$ given by
$$L_j^2(\xi,\eta)= (-\xi, -s_j\eta)~,$$
where $s_j$ is the slope of the line $\ell_j$ passing through the points $(a_{j+1},b_{j+1})$ and $(a_j,b_j)$. In particular, it follows from our hypothesis that  $0<s_j<1$. This is the main difference with respect to the setting in \cite{DG12} justifying the different asymmetry in the exponents of our result. 
For a square $S\in {\mathbf{S}}$ let $R_{S,j}$ denote the rectangle
 $$R_{S,j}:=(a_j,b_j)+L_j^2(S)~.$$
We define ${\mathbf{S}_j}$ to be the subcollection of squares $S\in\mathbf{S}$ such that $\alpha R_{S,j}\cap T_j \neq \emptyset$ where $\alpha <1$ is a fixed parameter to be chosen sufficiently close to $1$. In particular, the rectangles $\alpha R_{S,j}$ for $S\in {\mathbf{S}_j}$ cover $T_j\setminus \ell_j$. Moreover, choosing the constant $C_0$ above sufficiently large ensures that each $R_{S,j}$ lies inside the epigraph of our polygonal curve. We define
$$\mathbf{R}_j:=\lbrace R_{S,j}: S\in\mathbf{S}_j\rbrace, \qquad \mathbf{R}:= \bigcup_{j\ge 1}\mathbf{R}_j~.$$
For a rectangle $R\in\mathbf{R}$ let $J_R^1$ and $J_R^2$ denote the edges of $R$ parallel, respectively, to the horizontal and vertical axes. We define
$$J_j^1:= \bigcup_{R \in \mathbf{R}_j}J_R^1,\quad J_j^2:= \bigcup_{R \in \mathbf{R}_j}J_R^2~, \quad J_j^3:= \bigcup_{R \in \mathbf{R}_j}(-J_R^1-J_R^2)~,$$
and set $I^i_j:=\tfrac{1}{\alpha}J_j^i$, $i=1,2,3$. Choosing $\alpha$ sufficiently close to $1$ guarantees that the intervals in each collection $\lbrace I_j^i\rbrace_{j\in\mathbb{N}_0}$ have bounded overlap. Let $\varphi_j^i$, $i=1,2,3$, be fixed functions adapted to and supported on $I_j^i$. For later use we introduce the notation ${f}_j:=(\widehat{f}\varphi_j^1)^\vee$, ${g}_j:=(\widehat{g}\varphi_j^2)^\vee$, ${h}_j:=(\widehat{h}\varphi_j^3)^\vee$. 

Recall that $m$ in the statement of Theorem \ref{thm_polygonal} is the indicator function of the epigraph of our polygonal curve. Let $\Gamma_{\mathrm{par}}$ denote the staircase set associated with the sequences $\lbrace a_j\rbrace_{j\in\mathbb{N}_0}$ $\lbrace b_j\rbrace_{j\in\mathbb{N}_0}$. For a subset $\mathrm{C}\subseteq \mathbb{R}^d$, $d\in\mathbb{N}$, let $\psi_\mathrm{C}$ denote a fixed function adapted to and supported on $\mathrm{C}$.
By invoking, for example, \cite[Proposition~8.1]{DG12} we have that up to a set of measure zero a decomposition like the following exists
$$ m= \psi_{\Gamma_\mathrm{par}}+\sum_{R\in\mathbf{R}}\psi_R.$$
Boundedness in the desired range of exponents \eqref{range_beyond_L2} for the bilinear multiplier operator associated with $\psi_{\Gamma_\mathrm{par}}$ follows from Theorem \ref{thm_general_paraproduct}. We therefore focus on the second term.

Let $\mathbf{Q}$ be the collection of dyadic cubes $Q\subset \mathbb{R}^3$ with centers in $2^{j-10}\mathbb{Z}^3$, of side-length $2^j$ for some integer $j$, and satisfying the Whitney condition 
\begin{align*}
    C_0\, Q\cap\lbrace (\xi,\xi):\, \xi\in\mathbb{R} \rbrace & = \emptyset \\
    10 C_0\, Q\cap\lbrace (\xi,\xi):\, \xi\in\mathbb{R} \rbrace & \neq \emptyset~.
\end{align*}
Let $\lbrace L_j^3\rbrace_{j\in\mathbb{N}_0}$ be the family of linear transformations on $\mathbb{R}^3$ given by
$$L_j^3(\xi,\eta,\theta)=(-\xi,-s_j\eta, (1+s_j)\theta)~,$$
and define the collection of cubes $\mathbf{Q}_j:=\lbrace L_j^3(Q):\, Q\in\mathbf{Q}\rbrace$.
For a rectangle $R\in\mathbf{R}$ let $S$ be the square in $\mathbf{S}_j$, with $S=I\times J$, such that $R=L^2_j(S)+(a_j,b_j)$. We have $L^2_j(S)=(-I)\times (-s_j J)$ and if $(\xi,\eta)\in L^2_j(S)$ then $-\xi-\eta \in K_R:= I+s_jJ$. Let $\varphi_{K_R}$ be 
such that $\mathbf{1}_{K_R}\leq \varphi_{K_R}\le \mathbf{1}_{\tfrac{1}{\alpha}K_R}$. In particular, one can write $\varphi_{K_R}$ as
$$\varphi_{K_R}= \sum_{\omega_3: (R-(a_j,b_j))\times\omega_3 \in \mathbf{Q}_j} \psi_{\omega_3}$$
for some Fourier multipliers $\psi_{\omega_3}$ supported on and adapted to $\omega_3$. 

For an interval $\omega\subset\mathbb{R}$ let $\pi_\omega$ be the multiplier operator $\widehat{\pi_\omega f}:=\psi_{\omega}\widehat{f}$. For $a\in\mathbb{R}$ define $\widehat{M_a f}(\xi):= \widehat{f}(\xi+a)$. With this notation, and again following \cite{DG12}, the adjoint form that we want to bound can be rewritten as
\begin{align*}
    &\int_{\xi+\eta+\theta=0} \widehat{f}(\xi)\widehat{g}(\eta)\widehat{h}(\theta) \sum_{R\in\mathbf{R}_j} \psi_R(\xi,\eta) \dd\xi\, \dd\eta \, \dd\theta\\
    &=\int_{\xi+\eta+\theta=0} \widehat{f}_j(\xi)\widehat{g}_j(\eta)\widehat{h}_j(\theta) \sum_{R\in\mathbf{R}_j} \psi_R(\xi,\eta) \varphi_{K_R}(\theta- (-a_j-b_j) ) \dd\xi\, \dd\eta \, \dd\theta\\
    &= \sum_{Q=\omega_1\times \omega_2\times \omega_3 \in \mathbf{Q}_j} \int \pi_{\omega_1}(M_{a_j} f_j) \pi_{\omega_2}(M_{b_j}g_j) \pi_{\omega_3}(M_{-a_j-b_j} h_j)~.
\end{align*}
We define $\mathbf{P}_j$ to be the collection of multi-tiles $\vec{P}=(P_1,P_2,P_3)$ identified by a dyadic interval $I_{\vec{P}}$ in space and by cubes $\omega_1\times \omega_2 \times \omega_3 \in \mathbf{Q}_j$ in frequency. Define $j_{\vec{P}}$ to be the number such that $2^{-2^{C_0}j_{\vec{P}}}=|I_{\vec{P}}|$. Let $\rho$ denote a fixed positive function with frequency support in the interval $[-4^{-2^{C_0}}, 4^{-2^{C_0}}]$, normalized in $L^1$, and with sufficiently good decay in space, and define $\rho_j(x):=2^{-2^{C_0}j}\rho(2^{-2^{C_0}j}x)$. For an interval $I\subset \mathbb{R}$ define $\chi_{I,j}:=\mathbf{1}_{I}\ast \rho_j$. Hence, with this notation, we have the partition of unity at scale $j_0$
$$1=\sum_{\vec{P}\in\mathbf{P}_j : j_{\vec{P}}=j_0} \chi_{I_{\vec{p},j_{\vec{P}}}}~.$$
Relying on this, we can discretize in space the adjoint form obtaining the model sum 
\begin{align*}
    \sum_j \bigg| \sum_{\vec{P}\in\mathbf{P}_j} \int \chi_{I_{\vec{P},j_{\vec{P}}}}\pi_{\omega_1}(M_{a_j} f_j) \pi_{\omega_2}(M_{b_j}g_j) \pi_{\omega_3}(M_{-a_j-b_j} h_j)  \bigg| ~.
\end{align*}
Boundedness of this model sum form follows from \cite[Theorem~2.3]{DG12} by reversing the role of $\xi$ and $\eta$.
By duality, one recovers the boundedness of the multiplier operator in the range of exponents \eqref{range_beyond_L2}.
\qed

\section*{Acknowledgments}
The author would like to thank Christoph Thiele for introducing her to this topic and for inspiring discussions on the subject. She sincerely thanks Lars Becker for valuable discussions and comments on an earlier draft of this note, 
and Fred Yu-Hsiang Lin and Ioannis Parissis for helpful opinions and feedback.
She also thanks Gennady Uraltsev for interesting discussions on a related topic.
Part of this work has been carried out during the author's research stay at the Hausdorff Research Institute for Mathematics in Bonn, during the trimester program “Boolean Analysis in Computer Science”.


\begin{bibsection}
 \begin{biblist}
 \bib{Ba21}{article}{
author={Bakas, Odysseas},
   title={Sharp asymptotic estimates for a class of {L}ittlewood-{P}aley
   operators},
   journal={Studia Math.},
   volume={260},
   date={2021},
   number={2},
   pages={195--206},
   issn={0039-3223},
}
\bib{Bou89}{incollection}{
    AUTHOR = {Bourgain, Jean},
     TITLE = {On the behavior of the constant in the {L}ittlewood-{P}aley
              inequality},
 BOOKTITLE = {Geometric aspects of functional analysis (1987--88)},
    SERIES = {Lecture Notes in Math.},
    VOLUME = {1376},
     PAGES = {202--208},
 PUBLISHER = {Springer, Berlin},
      YEAR = {1989},
      ISBN = {3-540-51303-5},
}
\bib{CHL24}{article}{
  title={A sharp {H}örmander condition for bilinear {F}ourier multipliers with {L}ipschitz singularities},
  author={Chen, Jiao}, 
  AUTHOR = {Hsu, Martin}, 
  AUTHOR = {Lin, Fred Yu-Hsiang},
  journal={arXiv preprint arXiv:2403.04721},
  year={2024}
}
\bib{DG12}{article}{
    AUTHOR = {Demeter, Ciprian}, 
    AUTHOR = {Gautam, S. Zubin},
     TITLE = {Bilinear {F}ourier restriction theorems},
   JOURNAL = {J. Fourier Anal. Appl.},
    VOLUME = {18},
      YEAR = {2012},
    NUMBER = {6},
     PAGES = {1265--1290},
      ISSN = {1069-5869,1531-5851},
}
\bib{DPT16}{article}{
    AUTHOR = {Di Plinio, Francesco},
    AUTHOR = {Thiele, Christoph},
     TITLE = {Endpoint bounds for the bilinear {H}ilbert transform},
   JOURNAL = {Trans. Amer. Math. Soc.},
    VOLUME = {368},
      YEAR = {2016},
    NUMBER = {6},
     PAGES = {3931--3972},
      ISSN = {0002-9947,1088-6850},
}
\bib{DG07}{article}{
    AUTHOR = {Diestel, Geoff},
    AUTHOR = {Grafakos, Loukas},
     TITLE = {Unboundedness of the ball bilinear multiplier operator},
   JOURNAL = {Nagoya Math. J.},
    VOLUME = {185},
      YEAR = {2007},
     PAGES = {151--159},
      ISSN = {0027-7630,2152-6842},
}
\bib{G12}{article}{
    AUTHOR = {Gautam, S. Zubin},
     TITLE = {On curvature and the bilinear multiplier problem},
   JOURNAL = {Rev. Mat. Iberoam.},
    VOLUME = {28},
      YEAR = {2012},
    NUMBER = {2},
     PAGES = {351--369},
      ISSN = {0213-2230,2235-0616},
}
\bib{GN00}{article}{
    AUTHOR = {Gilbert, John E.}, 
    AUTHOR = {Nahmod, Andrea R.},
     TITLE = {Boundedness of bilinear operators with nonsmooth symbols},
   JOURNAL = {Math. Res. Lett.},
    VOLUME = {7},
      YEAR = {2000},
    NUMBER = {5-6},
     PAGES = {767--778},
      ISSN = {1073-2780},
}
\bib{GL04}{article}{
    AUTHOR = {Grafakos, Loukas},
    AUTHOR = {Li, Xiaochun},
     TITLE = {Uniform bounds for the bilinear {H}ilbert transforms. {I}},
   JOURNAL = {Ann. of Math. (2)},
    VOLUME = {159},
      YEAR = {2004},
    NUMBER = {3},
     PAGES = {889--933},
      ISSN = {0003-486X,1939-8980},
}
\bib{GL06}{article}{
    AUTHOR = {Grafakos, Loukas},
    AUTHOR = {Li, Xiaochun},
     TITLE = {The disc as a bilinear multiplier},
   JOURNAL = {Amer. J. Math.},
    VOLUME = {128},
      YEAR = {2006},
    NUMBER = {1},
     PAGES = {91--119},
      ISSN = {0002-9327,1080-6377},
}
\bib{GR10}{article}{
    AUTHOR = {Grafakos, Loukas}, 
    AUTHOR = {Reguera Rodr\'iguez, Maria Carmen},
     TITLE = {The bilinear multiplier problem for strictly convex compact  sets},
   JOURNAL = {Forum Math.},
    VOLUME = {22},
      YEAR = {2010},
    NUMBER = {4},
     PAGES = {619--626},
      ISSN = {0933-7741,1435-5337},
}
\bib{HK89}{article}{
    AUTHOR = {Hare, Kathryn E.},
    AUTHOR = {Klemes, Ivo},
     TITLE = {Properties of {L}ittlewood-{P}aley sets},
   JOURNAL = {Math. Proc. Cambridge Philos. Soc.},
    VOLUME = {105},
      YEAR = {1989},
    NUMBER = {3},
     PAGES = {485--494},
      ISSN = {0305-0041,1469-8064},
}
\bib{LT97}{article}{
    AUTHOR = {Lacey, Michael}, 
    AUTHOR = {Thiele, Christoph},
     TITLE = {{$L^p$} estimates for the bilinear {H}ilbert transform},
   JOURNAL = {Proc. Nat. Acad. Sci. U.S.A.},
    VOLUME = {94},
      YEAR = {1997},
    NUMBER = {1},
     PAGES = {33--35},
      ISSN = {0027-8424},
}
\bib{LT97a}{article}{
    AUTHOR = {Lacey, Michael}, 
    AUTHOR = {Thiele, Christoph},
     TITLE = {{$L^p$} estimates on the bilinear {H}ilbert transform for {$2<p<\infty$}},
   JOURNAL = {Ann. of Math. (2)},
    VOLUME = {146},
      YEAR = {1997},
    NUMBER = {3},
     PAGES = {693--724},
      ISSN = {0003-486X,1939-8980},
}
\bib{LT98}{article}{
    AUTHOR = {Lacey, Michael}, 
    AUTHOR = {Thiele, Christoph},
     TITLE = {On {C}alder\'on's conjecture for the bilinear {H}ilbert  transform},
   JOURNAL = {Proc. Nat. Acad. Sci. USA},
    VOLUME = {95},
      YEAR = {1998},
    NUMBER = {9},
     PAGES = {4828--4830},
      ISSN = {0027-8424,1091-6490},
}
\bib{LT99}{article}{
    AUTHOR = {Lacey, Michael},
    AUTHOR = {Thiele, Christoph},
     TITLE = {On {C}alder\'on's conjecture},
   JOURNAL = {Ann. of Math. (2)},
    VOLUME = {149},
      YEAR = {1999},
    NUMBER = {2},
     PAGES = {475--496},
      ISSN = {0003-486X,1939-8980},
}
\bib{L06}{article}{
    AUTHOR = {Li, Xiaochun},
     TITLE = {Uniform bounds for the bilinear {H}ilbert transforms. {II}},
   JOURNAL = {Rev. Mat. Iberoam.},
    VOLUME = {22},
      YEAR = {2006},
    NUMBER = {3},
     PAGES = {1069--1126},
      ISSN = {0213-2230,2235-0616},
}
\bib{L08}{article}{
    AUTHOR = {Li, Xiaochun},
     TITLE = {Uniform estimates for some paraproducts},
   JOURNAL = {New York J. Math.},
    VOLUME = {14},
      YEAR = {2008},
     PAGES = {145--192},
      ISSN = {1076-9803},
}
\bib{MS13}{book}{
    AUTHOR = {Muscalu, Camil},
    AUTHOR = {Schlag, Wilhelm},
     TITLE = {Classical and multilinear harmonic analysis. {V}ol. {II}},
    SERIES = {Cambridge Studies in Advanced Mathematics},
    VOLUME = {138},
 PUBLISHER = {Cambridge University Press, Cambridge},
      YEAR = {2013},
     PAGES = {xvi+324},
      ISBN = {978-1-107-03182-1},
}
\bib{MTT02}{incollection}{
    AUTHOR = {Muscalu, Camil},
    AUTHOR = {Tao, Terence}, 
    AUTHOR = {Thiele, Christoph},
     TITLE = {Uniform estimates on multi-linear operators with modulation symmetry},
      NOTE = {Dedicated to the memory of Tom Wolff},
   JOURNAL = {J. Anal. Math.},
    VOLUME = {88},
      YEAR = {2002},
     PAGES = {255--309},
      ISSN = {0021-7670,1565-8538},
}
\bib{M00}{book}{
    AUTHOR = {Muscalu, Florin Camil},
     TITLE = {L(p) estimates for multilinear operators given by singular symbols},
      NOTE = {Thesis (Ph.D.)--Brown University},
 PUBLISHER = {ProQuest LLC, Ann Arbor, MI},
      YEAR = {2000},
     PAGES = {145},
      ISBN = {978-0599-93764-2},
}
\bib{ST23}{article}{
    AUTHOR = {Saari, Olli},
    AUTHOR= {Thiele, Christoph},
     TITLE = {Paraproducts for bilinear multipliers associated with convex sets},
   JOURNAL = {Math. Ann.},
    VOLUME = {385},
      YEAR = {2023},
    NUMBER = {3-4},
     PAGES = {2013--2036},
      ISSN = {0025-5831,1432-1807},
}
\bib{T02}{article}{
    AUTHOR = {Thiele, Christoph},
     TITLE = {A uniform estimate},
   JOURNAL = {Ann. of Math. (2)},
    VOLUME = {156},
      YEAR = {2002},
    NUMBER = {2},
     PAGES = {519--563},
      ISSN = {0003-486X,1939-8980},
}
\bib{UW22}{article}{
  title={The full range of uniform bounds for the bilinear {H}ilbert transform},
  author={Uraltsev, Gennady},
  AUTHOR ={Warchalski, Micha{\l}},
  journal={arXiv preprint arXiv:2205.09851},
  year={2022}
}
 \end{biblist}
 \end{bibsection}
\end{document}